\documentclass[11pt]{article}

\usepackage[a4paper,margin=30mm]{geometry}
\usepackage{amsmath,amssymb,amsthm,mathtools}
\usepackage{microtype}
\usepackage{booktabs}
\usepackage{enumitem}
\usepackage[hidelinks]{hyperref}
\hypersetup{
  pdftitle={Erratum to Higher order scrambled digital nets achieve the optimal rate of the root mean square error for smooth integrands},
  pdfauthor={Josef Dick}
}
\usepackage{xcolor}

\newtheorem{theorem}{Theorem}[section]
\newtheorem{lemma}[theorem]{Lemma}
\newtheorem{proposition}[theorem]{Proposition}
\newtheorem{corollary}[theorem]{Corollary}
\newtheorem{remark}[theorem]{Remark}

\newtheorem{example}[theorem]{Example}

\newcommand{\N}{\mathbb N}
\newcommand{\Nzero}{\mathbb N_0}

\newcommand{\D}{\mathcal D}
\newcommand{\E}{\mathcal E}
\newcommand{\wal}{\operatorname{wal}}
\newcommand{\supp}{\operatorname{supp}}
\newcommand{\Var}{\operatorname{Var}}
\newcommand{\Vol}{\operatorname{Vol}}
\newcommand{\one}{\mathbf 1}
\newcommand{\abs}[1]{\lvert #1\rvert}
\newcommand{\norm}[1]{\lVert #1\rVert}
\newcommand{\set}[1]{\{#1\}}
\newcommand{\floor}[1]{\lfloor #1\rfloor}
\newcommand{\dd}{\,\mathrm d}
\newcommand{\eps}{\varepsilon}

\title{Erratum to ``Higher order scrambled digital nets achieve the optimal rate of the root mean square error for smooth integrands''}
\author{Josef Dick}
\date{1 August 2026}

\begin{document}
\maketitle

\begin{abstract}
This note corrects several statements and proof steps in J.~Dick,
\emph{Annals of Statistics} \textbf{39} (2011), 1372--1398,
DOI: 10.1214/11-AOS880, arXiv:1007.0842.
The principal smooth-function result remains valid: an order-$d$
nested-uniformly scrambled digital $(t,m,ds)$-net has root mean square
integration error $O(N^{-\min(\alpha,d)-1/2+\eps})$ for every $\eps>0$ for
functions in the unanchored Sobolev space with square-integrable mixed
partial derivatives up to order $\alpha$ in each variable.  The finite-difference variation introduced in the paper neither coincides with the displayed Sobolev norm nor directly
controls the finite differences used in Appendix~B; the claimed theorem for
that variation is therefore withdrawn here.  In addition, the variance bound
in Theorem~10 is missing a square, and the proof for $d>\alpha$ does not prove
the logarithmic power printed in the theorem.  We give a direct Sobolev proof, and a corrected logarithmic factor.
\end{abstract}

\section{Scope of the correction}

The following parts of the original paper are replaced or amended:
\begin{enumerate}[label=(\roman*),leftmargin=2.2em]
\item the definition of order-$d$ scrambling and Proposition~5;
\item the covariance statement used in Lemmas~6 and~7;
\item equation~(3.2) and the limiting justification in Lemma~7;
\item the finite-difference variation on pp.~1386--1387 and its use in Appendix~B;
\item Lemma~9 and its proof in Appendix~B, replaced by a direct Sobolev estimate;
\item Theorem~10 and the logarithmic factor stated in Section~4;
\item the typographical and indexing corrections collected in Section~\ref{sec:line-corrections}.
\end{enumerate}
The numerical experiments and the optimal algebraic exponent
$\min(\alpha,d)+1/2$ for the mixed Sobolev class are unchanged.  No claim
about the original finite-difference variation is made in this erratum.

Throughout, $b$ is a prime, $s,d,\alpha,m\in\N$, $0\le t\le m$,
$N=b^m$, and
\[
  M:=m-t,\qquad \rho:=\min(\alpha,d).
\]
All base-$b$ expansions of points in $[0,1)$ are chosen to be the expansion that
is not eventually equal to $b-1$.

\section{Digit interlacing and scrambling}

\subsection{The order-\texorpdfstring{$d$}{d} scrambling is defined before interlacing}

The original paper writes $\D_d(\Pi(\D_d^{-1}(x)))$ for an arbitrary
$x\in[0,1)^s$.  This is not well-defined because Lemma~11 correctly states that
$\D_d$ is not surjective in general, but is surjective if a set of Lebesgue measure $0$ is removed.  No inverse is required by the algorithm.

For $z=(z_1,\ldots,z_{ds})\in[0,1)^{ds}$, write
$z_j=\sum_{a\ge1}z_{j,a}b^{-a}$ and define
\[
 \D_d(z)_i
 :=\sum_{a=1}^{\infty}\sum_{r=1}^{d}
     z_{(i-1)d+r,a}\,b^{-r-(a-1)d},
 \qquad 1\le i\le s.
\]
Let $\Pi$ be Owen's nested uniform scrambling, applied independently to the
$ds$ coordinates.  The order-$d$ scrambled point is
\begin{equation}\label{eq:correct-order-d-scramble}
  Y:=\D_d(\Pi(z)).
\end{equation}

\begin{proposition}[Replacement for Proposition 5]\label{prop:uniform}
For every deterministic $z\in[0,1)^{ds}$, the random point
$Y=\D_d(\Pi(z))$ is uniformly distributed on $[0,1)^s$.
\end{proposition}

\begin{proof}
Each coordinate of $\Pi(z)$ is uniformly distributed on $[0,1)$, and the
$ds$ coordinates are independent.  More explicitly, fix $L\ge1$ and prescribe
the first $dL$ digits of every output coordinate.  By the definition of
interlacing, this is exactly the same as prescribing the first $L$ digits of
each of the $ds$ scrambled input coordinates.  Those $dsL$ digits are
independent and uniform, so every such output digit pattern has probability
$b^{-dsL}$.

It follows that every half-open $s$-dimensional base-$b$ box whose side length
is $b^{-dL}$ has probability equal to its volume.  To pass from these grid
boxes to an arbitrary rectangle
$[0,x_1)\times\cdots\times[0,x_s)$, round each $x_i$ down and up to the nearest
multiple of $b^{-dL}$.  The desired probability is squeezed between the
volumes of the two resulting grid rectangles.  As $L\to\infty$, both volumes
converge to $x_1\cdots x_s$.  Thus the joint distribution function of $Y$ is
the distribution function of the uniform measure on $[0,1)^s$.
\end{proof}

\subsection{A covariance statement without an inverse interlacing map}

For $\ell\in\Nzero$, put
\[
 B_0:=\set{0},\qquad
 B_\ell:=\set{k\in\Nzero:b^{\ell-1}\le k<b^\ell}\quad(\ell\ge1).
\]
For $x,x'\in[0,1)$ define
\[
 \chi_0(x,x'):=1,
\]
and, for $\ell\ge1$,
\begin{equation}\label{eq:chi}
 \chi_\ell(x,x')
 :=\frac{b\,\one_{\{\floor{b^\ell x}=\floor{b^\ell x'}\}}
       -\one_{\{\floor{b^{\ell-1}x}=\floor{b^{\ell-1}x'}\}}}{b-1}.
\end{equation}
For vectors, $B_{\boldsymbol\ell}:=\prod_{j=1}^{ds}B_{\ell_j}$.

\begin{proposition}[Covariance under nested uniform scrambling]\label{prop:covariance}
Let $z,z'\in[0,1)^{ds}$ be deterministic and let the same nested uniform
scrambling $\Pi$ be applied to both points.  For
$\boldsymbol k,\boldsymbol k'\in\Nzero^{ds}$,
\[
 \mathbb E\!\left[
 \wal_{\boldsymbol k}(\Pi(z))
 \overline{\wal_{\boldsymbol k'}(\Pi(z'))}
 \right]=0
 \quad\text{if }\boldsymbol k\ne\boldsymbol k'.
\]
If $\boldsymbol k=\boldsymbol k'\in B_{\boldsymbol\ell}$, then
\begin{equation}\label{eq:block-covariance}
 \mathbb E\!\left[
 \wal_{\boldsymbol k}(\Pi(z))
 \overline{\wal_{\boldsymbol k}(\Pi(z'))}
 \right]
 =\prod_{j=1}^{ds}\chi_{\ell_j}(z_j,z'_j).
\end{equation}
In particular, the diagonal covariance depends on the frequency only through
its block index $\boldsymbol\ell$.
\end{proposition}

\begin{proof}
We first prove the corresponding one-dimensional statement and then use
independence of the scrambling permutations between coordinates.

Let \(x,x'\in[0,1)\), and write their base-\(b\) expansions as
\[
 x=\sum_{r=1}^\infty x_r b^{-r},
 \qquad
 x'=\sum_{r=1}^\infty x'_r b^{-r},
 \qquad
 x_r,x'_r\in\{0,\ldots,b-1\}.
\]
For \(k,k'\in\Nzero\), write
\[
 k=\sum_{r=0}^\infty \kappa_r b^r,
 \qquad
 k'=\sum_{r=0}^\infty \kappa'_r b^r,
 \qquad
 \kappa_r,\kappa'_r\in\{0,\ldots,b-1\},
\]
where only finitely many digits are nonzero.  Let
\(\omega_b=\exp(2\pi i/b)\).  Then
\[
 \wal_k(x)
 =
 \omega_b^{\sum_{r\geq 1}\kappa_{r-1}x_r}.
\]

Under nested uniform scrambling, for every finite digit string
\((a_1,\ldots,a_m)\) there is an independent uniformly distributed
permutation
\[
 \pi_{a_1,\ldots,a_m}
\]
of \(\{0,\ldots,b-1\}\).  Thus, if
\(\widetilde x=\Pi(x)\) and \(\widetilde x'=\Pi(x')\), then
\[
 \widetilde x_r
 =
 \pi_{x_1,\ldots,x_{r-1}}(x_r),
 \qquad
 \widetilde x'_r
 =
 \pi_{x'_1,\ldots,x'_{r-1}}(x'_r).
\]
The same family of permutations is used for both \(x\) and \(x'\).

We shall repeatedly use the following elementary fact.  If \(\pi\) is a
uniform random permutation of \(\{0,\ldots,b-1\}\), then for any
\(a\in\{0,\ldots,b-1\}\) and any \(u\in\{0,\ldots,b-1\}\),
\[
 \mathbb E\!\left[\omega_b^{u\pi(a)}\right]
 =
 \begin{cases}
  1, & u=0,\\
  0, & u\neq 0.
 \end{cases}
\]
Moreover, if \(a\neq a'\), then \((\pi(a),\pi(a'))\) is uniformly
distributed over all ordered pairs of distinct digits, and hence
\begin{equation}\label{eq:perm-pair-expectation}
 \mathbb E\!\left[
   \omega_b^{u\pi(a)-v\pi(a')}
 \right]
 =
 \begin{cases}
  1, & u=v=0,\\[1mm]
  -\dfrac{1}{b-1}, & u=v\neq 0,\\[2mm]
  0, & u\neq v.
 \end{cases}
\end{equation}
Indeed,
\[
 \frac{1}{b(b-1)}
 \sum_{\substack{c,d=0\\c\neq d}}^{b-1}
 \omega_b^{uc-vd}
\]
is the required expectation, and the displayed values follow from
\[
 \sum_{c=0}^{b-1}\omega_b^{mc}
 =
 \begin{cases}
  b, & m=0\pmod b,\\
  0, & m\neq 0\pmod b.
 \end{cases}
\]

We now consider
\[
 \mathbb E\!\left[
   \wal_k(\Pi(x))
   \overline{\wal_{k'}(\Pi(x'))}
 \right].
\]
Since
\[
 \wal_k(\Pi(x))
 \overline{\wal_{k'}(\Pi(x'))}
 =
 \prod_{r\geq 1}
 \omega_b^{
   \kappa_{r-1}\widetilde x_r
   -
   \kappa'_{r-1}\widetilde x'_r
 },
\]
suppose first that \(k\neq k'\).  Then there is a digit position \(r\)
such that
\[
 \kappa_{r-1}\neq \kappa'_{r-1}.
\]
Condition on all scrambling permutations except the permutation or
permutations used to generate the \(r\)-th scrambled digits.  There are
three possibilities.

If
\[
 (x_1,\ldots,x_{r-1})=(x'_1,\ldots,x'_{r-1})
 \quad\text{and}\quad
 x_r=x'_r,
\]
then the same permutation is applied to the same digit and the
\(r\)-th factor is
\[
 \omega_b^{(\kappa_{r-1}-\kappa'_{r-1})\pi(x_r)}.
\]
Its conditional expectation is zero because
\(\kappa_{r-1}-\kappa'_{r-1}\not\equiv 0\pmod b\).

If
\[
 (x_1,\ldots,x_{r-1})=(x'_1,\ldots,x'_{r-1})
 \quad\text{but}\quad
 x_r\neq x'_r,
\]
then the same random permutation is applied to two different digits.
By \eqref{eq:perm-pair-expectation}, the corresponding conditional
expectation is zero because
\(\kappa_{r-1}\neq\kappa'_{r-1}\).

Finally, if
\[
 (x_1,\ldots,x_{r-1})\neq(x'_1,\ldots,x'_{r-1}),
\]
then the two \(r\)-th scrambled digits are generated by independent
random permutations.  Hence the corresponding conditional expectation
factors as
\[
 \mathbb E\!\left[
   \omega_b^{\kappa_{r-1}\widetilde x_r}
 \right]
 \mathbb E\!\left[
   \omega_b^{-\kappa'_{r-1}\widetilde x'_r}
 \right].
\]
Since
\(\kappa_{r-1}\neq\kappa'_{r-1}\), at least one of these two digits is
nonzero, so at least one factor is zero.  Thus in every case
\[
 \mathbb E\!\left[
   \wal_k(\Pi(x))
   \overline{\wal_{k'}(\Pi(x'))}
 \right]
 =0
 \qquad\text{whenever }k\neq k'.
\]

It remains to consider the diagonal case.  Let \(k=k'\in B_\ell\).
For \(\ell\geq 1\), this means that
\[
 \kappa_{\ell-1}\neq 0,
 \qquad
 \kappa_r=0
 \quad\text{for all }r\geq \ell.
\]
Let \(q\) be the number of initial base-\(b\) digits shared by \(x\)
and \(x'\), that is,
\[
 x_1=x'_1,\ldots,x_q=x'_q,
\]
and either \(x_{q+1}\neq x'_{q+1}\) or \(q=\infty\) if \(x=x'\).

If \(q\geq\ell\), then the first \(\ell\) scrambled digits of \(x\)
and \(x'\) agree, so
\[
 \wal_k(\Pi(x))
 \overline{\wal_k(\Pi(x'))}
 =1
\]
and therefore the expectation equals \(1\).

If \(q=\ell-1\), then the two points agree through digit
\(\ell-1\) and differ at digit \(\ell\).  The factors corresponding to
digits \(1,\ldots,\ell-1\) cancel.  At digit \(\ell\), the same random
permutation is applied to two distinct digits, and since
\(\kappa_{\ell-1}\neq 0\), \eqref{eq:perm-pair-expectation} gives
\[
 \mathbb E\!\left[
 \omega_b^{\kappa_{\ell-1}
   \left(\widetilde x_\ell-\widetilde x'_\ell\right)}
 \right]
 =
 -\frac{1}{b-1}.
\]
Hence in this case the covariance equals \(-1/(b-1)\).

If \(q\leq\ell-2\), then the prefixes of length \(\ell-1\) of the two
points are already different.  Therefore their \(\ell\)-th scrambled
digits are generated by independent random permutations.  Since
\(\kappa_{\ell-1}\neq0\),
\[
 \mathbb E\!\left[
   \omega_b^{\kappa_{\ell-1}\widetilde x_\ell}
 \right]
 =
 0,
\]
and consequently the covariance is zero.

Thus, for every \(k\in B_\ell\),
\[
 \mathbb E\!\left[
   \wal_k(\Pi(x))
   \overline{\wal_k(\Pi(x'))}
 \right]
 =
 \begin{cases}
  1,
  &
  \lfloor b^\ell x\rfloor
  =
  \lfloor b^\ell x'\rfloor,
  \\[1mm]
  -\dfrac{1}{b-1},
  &
  \lfloor b^{\ell-1}x\rfloor
  =
  \lfloor b^{\ell-1}x'\rfloor
  \text{ and }
  \lfloor b^\ell x\rfloor
  \neq
  \lfloor b^\ell x'\rfloor,
  \\[3mm]
  0,
  &
  \lfloor b^{\ell-1}x\rfloor
  \neq
  \lfloor b^{\ell-1}x'\rfloor.
 \end{cases}
\]
By the definition in \eqref{eq:chi}, this is precisely
\[
 \chi_\ell(x,x').
\]
For \(\ell=0\), we have \(k=0\) and
\(\wal_0\equiv1\), so the same conclusion holds with
\(\chi_0(x,x')=1\).  We have therefore proved the one-dimensional
identity
\begin{equation}\label{eq:one-dimensional-scramble-covariance}
 \mathbb E\!\left[
   \wal_k(\Pi(x))
   \overline{\wal_{k'}(\Pi(x'))}
 \right]
 =
 \begin{cases}
  0, & k\neq k',\\
  \chi_\ell(x,x'), & k=k'\in B_\ell.
 \end{cases}
\end{equation}

We now return to \(ds\) dimensions.  Write
\[
 z=(z_1,\ldots,z_{ds}),
 \qquad
 z'=(z'_1,\ldots,z'_{ds}),
\]
and
\[
 \boldsymbol k=(k_1,\ldots,k_{ds}),
 \qquad
 \boldsymbol k'=(k'_1,\ldots,k'_{ds}).
\]
Since multivariate Walsh functions factor coordinatewise,
\[
 \wal_{\boldsymbol k}(\Pi(z))
 \overline{\wal_{\boldsymbol k'}(\Pi(z'))}
 =
 \prod_{j=1}^{ds}
 \wal_{k_j}(\Pi_j(z_j))
 \overline{\wal_{k'_j}(\Pi_j(z'_j))}.
\]
The nested scrambling permutation families used in different
coordinates are independent.  Hence
\[
 \mathbb E\!\left[
 \wal_{\boldsymbol k}(\Pi(z))
 \overline{\wal_{\boldsymbol k'}(\Pi(z'))}
 \right]
 =
 \prod_{j=1}^{ds}
 \mathbb E\!\left[
 \wal_{k_j}(\Pi_j(z_j))
 \overline{\wal_{k'_j}(\Pi_j(z'_j))}
 \right].
\]

If \(\boldsymbol k\neq\boldsymbol k'\), then
\(k_j\neq k'_j\) for at least one coordinate \(j\), and the
corresponding factor vanishes by
\eqref{eq:one-dimensional-scramble-covariance}.  Therefore
\[
 \mathbb E\!\left[
 \wal_{\boldsymbol k}(\Pi(z))
 \overline{\wal_{\boldsymbol k'}(\Pi(z'))}
 \right]
 =0.
\]

Finally, if
\(\boldsymbol k=\boldsymbol k'\in B_{\boldsymbol\ell}\), then
\(k_j\in B_{\ell_j}\) for every \(j\), and
\eqref{eq:one-dimensional-scramble-covariance} gives
\[
 \mathbb E\!\left[
 \wal_{\boldsymbol k}(\Pi(z))
 \overline{\wal_{\boldsymbol k}(\Pi(z'))}
 \right]
 =
 \prod_{j=1}^{ds}
 \chi_{\ell_j}(z_j,z'_j).
\]
The right-hand side depends on \(\boldsymbol k\) only through its block
index \(\boldsymbol\ell\), which proves the final assertion.
\end{proof}

This proposition replaces the indexing in Lemma~6.  It also handles identical
points without a special ``largest common digit'' convention. 

\section{The variance decomposition}

Let $\E_d:\Nzero^{ds}\to\Nzero^s$ be the digit-interlacing bijection on
frequency indices.  If
$k_j=\sum_{a\ge0}\kappa_{j,a}b^a$, then
\[
 \E_d(\boldsymbol k)_i
 =\sum_{a\ge0}\sum_{r=1}^{d}
   \kappa_{(i-1)d+r,a}b^{r-1+ad}.
\]
The identity
\begin{equation}\label{eq:wal-interlace}
 \wal_{\E_d(\boldsymbol k)}(\D_d(z))=\wal_{\boldsymbol k}(z)
\end{equation}
holds for the canonical digit expansions.

For $f\in L_2([0,1]^s)$ define
\[
 \widehat f(\boldsymbol h)
 :=\int_{[0,1]^s}f(x)\overline{\wal_{\boldsymbol h}(x)}\dd x
\]
and
\begin{equation}\label{eq:sigma-def}
 \sigma_{d,\boldsymbol\ell,s}^2(f)
 :=\sum_{\boldsymbol k\in B_{\boldsymbol\ell}}
       \abs{\widehat f(\E_d(\boldsymbol k))}^2.
\end{equation}
For an underlying deterministic point set
$P=\set{z_0,\ldots,z_{N-1}}\subset[0,1)^{ds}$ put
\begin{equation}\label{eq:gain-def}
 \Gamma_{d,\boldsymbol\ell}(P)
 :=\frac1{N^2}\sum_{n,n'=0}^{N-1}
       \prod_{j=1}^{ds}\chi_{\ell_j}(z_{n,j},z_{n',j}).
\end{equation}
For any $\boldsymbol k\in B_{\boldsymbol\ell}$,
Proposition~\ref{prop:covariance} also gives
\begin{equation}\label{eq:gain-square}
 \Gamma_{d,\boldsymbol\ell}(P)
 =\mathbb E\left|\frac1N\sum_{n=0}^{N-1}
   \wal_{\boldsymbol k}(\Pi(z_n))\right|^2.
\end{equation}
Consequently $0\le\Gamma_{d,\boldsymbol\ell}(P)\le1$.

Equation~(3.2) in the original paper must read
\begin{equation}\label{eq:parseval-var}
 \Var(f)=\sum_{\boldsymbol h\in\Nzero^s\setminus\{\boldsymbol0\}}
          \abs{\widehat f(\boldsymbol h)}^2,
\end{equation}
not the sum over all frequencies.

\begin{lemma}[Corrected Lemma 7]\label{lem:variance-decomposition}
Let $f\in L_2([0,1]^s)$ and
\[
 \widehat I_N(f):=\frac1N\sum_{n=0}^{N-1}
 f\bigl(\D_d(\Pi(z_n))\bigr).
\]
Then $\widehat I_N(f)$ is unbiased and
\begin{equation}\label{eq:variance-decomposition}
 \Var[\widehat I_N(f)]
 =\sum_{\boldsymbol\ell\in\Nzero^{ds}\setminus\{\boldsymbol0\}}
    \sigma_{d,\boldsymbol\ell,s}^2(f)
    \Gamma_{d,\boldsymbol\ell}(P).
\end{equation}
\end{lemma}

\begin{proof}
Let
\[
  g:[0,1)^{ds}\to\mathbb C,
  \qquad
  g(\boldsymbol y):=f\bigl(\D_d(\boldsymbol y)\bigr).
\]
The digit interlacing map \(\D_d\) is measure preserving: apart from the
null set on which a number has two base-\(b\) expansions, it merely
rearranges the base-\(b\) digits of the \(ds\) input coordinates into the
digits of the \(s\) output coordinates.  Consequently,
\[
  \lVert g\rVert_{L_2([0,1)^{ds})}^2
  =\int_{[0,1)^{ds}}
       \bigl|f(\D_d(\boldsymbol y))\bigr|^2\,\mathrm d\boldsymbol y
  =\int_{[0,1)^s}|f(\boldsymbol x)|^2\,\mathrm d\boldsymbol x
  =\lVert f\rVert_{L_2([0,1)^s)}^2.
\]
In particular, \(g\in L_2([0,1)^{ds})\).

Let
\(\E_d:\Nzero^{ds}\to\Nzero^s\)
denote the corresponding digit interlacing map on Walsh frequencies.
The defining digit rearrangements give
\[
  \wal_{\E_d(\boldsymbol k)}\bigl(\D_d(\boldsymbol y)\bigr)
  =\wal_{\boldsymbol k}(\boldsymbol y)
  \qquad
  \text{for all \(\boldsymbol k\in\Nzero^{ds}\) and a.e. \(\boldsymbol y\).}
\]
Since \(\D_d\) is measure preserving, the Walsh coefficients of \(g\)
therefore satisfy
\[
\begin{aligned}
  \widehat g(\boldsymbol k)
  &=\int_{[0,1)^{ds}}
       f\bigl(\D_d(\boldsymbol y)\bigr)
       \overline{\wal_{\boldsymbol k}(\boldsymbol y)}
       \,\mathrm d\boldsymbol y \\
  &=\int_{[0,1)^{ds}}
       f\bigl(\D_d(\boldsymbol y)\bigr)
       \overline{
         \wal_{\E_d(\boldsymbol k)}
         \bigl(\D_d(\boldsymbol y)\bigr)
       }
       \,\mathrm d\boldsymbol y \\
  &=\int_{[0,1)^s}
       f(\boldsymbol x)
       \overline{\wal_{\E_d(\boldsymbol k)}(\boldsymbol x)}
       \,\mathrm d\boldsymbol x
   =\widehat f\bigl(\E_d(\boldsymbol k)\bigr).
\end{aligned}
\]
In particular,
\[
  \widehat g(\boldsymbol0)
  =\widehat f(\boldsymbol0)
  =I(f):=\int_{[0,1)^s}f(\boldsymbol x)\,\mathrm d\boldsymbol x.
\]

\paragraph{Unbiasedness.}
For each fixed deterministic point \(z_n\), nested uniform scrambling
makes \(\Pi(z_n)\) uniformly distributed on \([0,1)^{ds}\).  Indeed, in
each coordinate the permutation used at a given digit is uniform, and,
conditionally on all preceding scrambled digits, the next scrambled digit
is uniform on \(\{0,\ldots,b-1\}\).  The permutation families belonging to
different coordinates are independent.  Thus every finite collection of
leading digits of \(\Pi(z_n)\) has the uniform distribution, which proves
the asserted uniformity.

Because \(\D_d\) is measure preserving,
\(\D_d(\Pi(z_n))\) is consequently uniform on \([0,1)^s\).  Moreover,
\(f\in L_2([0,1)^s)\subset L_1([0,1)^s)\), since the domain has measure
one.  Hence
\[
  \mathbb E\!\left[
    f\bigl(\D_d(\Pi(z_n))\bigr)
  \right]
  =\int_{[0,1)^s}f(\boldsymbol x)\,\mathrm d\boldsymbol x
  =I(f).
\]
It follows by linearity of expectation that
\[
\begin{aligned}
  \mathbb E[\widehat I_N(f)]
  &=\frac1N\sum_{n=0}^{N-1}
      \mathbb E\!\left[
        f\bigl(\D_d(\Pi(z_n))\bigr)
      \right] \\
  &=I(f).
\end{aligned}
\]
Thus \(\widehat I_N(f)\) is unbiased.  Notice that independence of the
scrambled points is neither assumed nor needed; only the marginal
uniformity of each scrambled point is used here.

\paragraph{Finite Walsh projections.}
The Walsh system is a complete orthonormal basis of
\(L_2([0,1)^{ds})\).  The coefficient identity above therefore gives the
\(L_2\)-expansion
\[
  g(\boldsymbol y)
  =\sum_{\boldsymbol k\in\Nzero^{ds}}
     \widehat f\bigl(\E_d(\boldsymbol k)\bigr)
     \wal_{\boldsymbol k}(\boldsymbol y).
\]
We first work with finite Walsh projections, thereby avoiding any formal
interchange of an infinite series and expectation.  For \(L\geq1\), put
\[
  A_L
  :=\left\{
      \boldsymbol k\in\Nzero^{ds}:
      0\leq k_j<b^L \text{ for }j=1,\ldots,ds
    \right\}.
\]
With the usual convention
\(B_0=\{0\}\) and
\(B_\ell=\{b^{\ell-1},\ldots,b^\ell-1\}\) for \(\ell\geq1\), this is
equivalently
\[
  A_L
  =\bigcup_{\boldsymbol\ell\in\{0,\ldots,L\}^{ds}}
     B_{\boldsymbol\ell}.
\]
Define
\[
  g_L(\boldsymbol y)
  :=\sum_{\boldsymbol k\in A_L}
      \widehat f\bigl(\E_d(\boldsymbol k)\bigr)
      \wal_{\boldsymbol k}(\boldsymbol y).
\]
Then \(g_L\to g\) in \(L_2([0,1)^{ds})\).  Since the constant
coefficient is retained,
\[
  \int_{[0,1)^{ds}}g_L(\boldsymbol y)\,\mathrm d\boldsymbol y=I(f),
\]
and hence
\[
  g_L(\boldsymbol y)-I(f)
  =\sum_{\substack{\boldsymbol k\in A_L\\
                    \boldsymbol k\neq\boldsymbol0}}
      \widehat f\bigl(\E_d(\boldsymbol k)\bigr)
      \wal_{\boldsymbol k}(\boldsymbol y).
\]
Set
\[
  \widehat I_{N,L}
  :=\frac1N\sum_{n=0}^{N-1}g_L\bigl(\Pi(z_n)\bigr).
\]
All sums in the following computation are finite.  Expanding the squared
absolute value and taking expectations gives
\[
\begin{aligned}
  \mathbb E\!\left[
    \left|\widehat I_{N,L}-I(f)\right|^2
  \right]
  ={}&\frac1{N^2}\sum_{n,n'=0}^{N-1}
      \sum_{\substack{\boldsymbol k,\boldsymbol k'\in A_L\\
                       \boldsymbol k,\boldsymbol k'\neq\boldsymbol0}}
      \widehat f\bigl(\E_d(\boldsymbol k)\bigr)
      \overline{
        \widehat f\bigl(\E_d(\boldsymbol k')\bigr)
      } \\
  &\quad\times
      \mathbb E\!\left[
        \wal_{\boldsymbol k}\bigl(\Pi(z_n)\bigr)
        \overline{
          \wal_{\boldsymbol k'}\bigl(\Pi(z_{n'})\bigr)
        }
      \right].
\end{aligned}
\]
By Proposition~\ref{prop:covariance}, the expectation in the last display
vanishes whenever \(\boldsymbol k\neq\boldsymbol k'\).  Thus all
off-diagonal Walsh-frequency terms disappear, leaving
\[
\begin{aligned}
  \mathbb E\!\left[
    \left|\widehat I_{N,L}-I(f)\right|^2
  \right]
  ={}&\frac1{N^2}\sum_{n,n'=0}^{N-1}
      \sum_{\substack{\boldsymbol k\in A_L\\
                       \boldsymbol k\neq\boldsymbol0}}
      \left|
        \widehat f\bigl(\E_d(\boldsymbol k)\bigr)
      \right|^2 \\
  &\quad\times
      \mathbb E\!\left[
        \wal_{\boldsymbol k}\bigl(\Pi(z_n)\bigr)
        \overline{
          \wal_{\boldsymbol k}\bigl(\Pi(z_{n'})\bigr)
        }
      \right].
\end{aligned}
\]
If \(\boldsymbol k\in B_{\boldsymbol\ell}\), the diagonal part of that
proposition yields
\[
  \mathbb E\!\left[
    \wal_{\boldsymbol k}\bigl(\Pi(z_n)\bigr)
    \overline{
      \wal_{\boldsymbol k}\bigl(\Pi(z_{n'})\bigr)
    }
  \right]
  =\prod_{j=1}^{ds}
     \chi_{\ell_j}(z_{n,j},z_{n',j}).
\]
The right-hand side depends on \(\boldsymbol k\) only through its block
index \(\boldsymbol\ell\).  We may therefore group the preceding sum by
Walsh blocks to obtain
\[
\begin{aligned}
  \mathbb E\!\left[
    \left|\widehat I_{N,L}-I(f)\right|^2
  \right]
  ={}&
  \sum_{\substack{
          \boldsymbol\ell\in\{0,\ldots,L\}^{ds}\\
          \boldsymbol\ell\neq\boldsymbol0}}
  \left(
    \sum_{\boldsymbol k\in B_{\boldsymbol\ell}}
      \left|
        \widehat f\bigl(\E_d(\boldsymbol k)\bigr)
      \right|^2
  \right) \\
  &\quad\times
  \left(
    \frac1{N^2}\sum_{n,n'=0}^{N-1}
      \prod_{j=1}^{ds}
        \chi_{\ell_j}(z_{n,j},z_{n',j})
  \right).
\end{aligned}
\]
By the definitions of the block energy and the gain coefficient,
\[
  \sigma_{d,\boldsymbol\ell,s}^2(f)
  =\sum_{\boldsymbol k\in B_{\boldsymbol\ell}}
     \left|
       \widehat f\bigl(\E_d(\boldsymbol k)\bigr)
     \right|^2
\]
and
\[
  \Gamma_{d,\boldsymbol\ell}(P)
  =\frac1{N^2}\sum_{n,n'=0}^{N-1}
     \prod_{j=1}^{ds}
       \chi_{\ell_j}(z_{n,j},z_{n',j}).
\]
Consequently,
\[
  \mathbb E\!\left[
    \left|\widehat I_{N,L}-I(f)\right|^2
  \right]
  =\sum_{\substack{
       \boldsymbol\ell\in\{0,\ldots,L\}^{ds}\\
       \boldsymbol\ell\neq\boldsymbol0}}
     \sigma_{d,\boldsymbol\ell,s}^2(f)
     \Gamma_{d,\boldsymbol\ell}(P).
\]
The zero block is absent because
\(B_{\boldsymbol0}=\{\boldsymbol0\}\) contains only the constant Walsh
frequency, which disappears upon centering by \(I(f)\).

\paragraph{Passage to the \(L_2\)-limit.}
For each \(n\), the point \(\Pi(z_n)\) is uniform on \([0,1)^{ds}\).
It follows that
\[
\begin{aligned}
  \mathbb E\!\left[
    \left|
      g_L\bigl(\Pi(z_n)\bigr)-g\bigl(\Pi(z_n)\bigr)
    \right|^2
  \right]
  &=\int_{[0,1)^{ds}}
      |g_L(\boldsymbol y)-g(\boldsymbol y)|^2
      \,\mathrm d\boldsymbol y \\
  &=\lVert g_L-g\rVert_2^2.
\end{aligned}
\]
The triangle inequality in \(L_2(\Omega)\) now gives
\[
\begin{aligned}
  \left\|
    \widehat I_{N,L}-\widehat I_N(f)
  \right\|_{L_2(\Omega)}
  &\leq\frac1N\sum_{n=0}^{N-1}
    \left\|
      g_L\bigl(\Pi(z_n)\bigr)-g\bigl(\Pi(z_n)\bigr)
    \right\|_{L_2(\Omega)} \\
  &=\lVert g_L-g\rVert_2
    \longrightarrow0.
\end{aligned}
\]
Therefore
\[
  \widehat I_{N,L}-I(f)
  \longrightarrow
  \widehat I_N(f)-I(f)
  \qquad\text{in }L_2(\Omega),
\]
and hence
\[
  \mathbb E\!\left[
    \left|\widehat I_{N,L}-I(f)\right|^2
  \right]
  \longrightarrow
  \mathbb E\!\left[
    \left|\widehat I_N(f)-I(f)\right|^2
  \right]
  =\Var[\widehat I_N(f)],
\]
where the last equality uses the unbiasedness already proved.

It remains only to identify the limit of the finite sums.  For every
\(\boldsymbol\ell\) and any
\(\boldsymbol k\in B_{\boldsymbol\ell}\),
Proposition~\ref{prop:covariance} implies
\[
\begin{aligned}
  \Gamma_{d,\boldsymbol\ell}(P)
  &=\frac1{N^2}\sum_{n,n'=0}^{N-1}
    \mathbb E\!\left[
      \wal_{\boldsymbol k}\bigl(\Pi(z_n)\bigr)
      \overline{
        \wal_{\boldsymbol k}\bigl(\Pi(z_{n'})\bigr)
      }
    \right] \\
  &=\mathbb E\!\left[
      \left|
        \frac1N\sum_{n=0}^{N-1}
          \wal_{\boldsymbol k}\bigl(\Pi(z_n)\bigr)
      \right|^2
    \right].
\end{aligned}
\]
Thus
\(0\leq\Gamma_{d,\boldsymbol\ell}(P)\leq1\),
because every Walsh function has modulus one.  Hence the finite sums above
increase, as \(L\to\infty\), to the full series.  This series is finite:
by Parseval's identity and the coefficient identity for \(g\),
\[
\begin{aligned}
  \sum_{\boldsymbol\ell\in\Nzero^{ds}}
    \sigma_{d,\boldsymbol\ell,s}^2(f)
    \Gamma_{d,\boldsymbol\ell}(P)
  &\leq
  \sum_{\boldsymbol\ell\in\Nzero^{ds}}
    \sigma_{d,\boldsymbol\ell,s}^2(f) \\
  &=\sum_{\boldsymbol k\in\Nzero^{ds}}
      \left|
        \widehat f\bigl(\E_d(\boldsymbol k)\bigr)
      \right|^2 \\
  &=\lVert g\rVert_2^2
   =\lVert f\rVert_2^2
   <\infty.
\end{aligned}
\]
Letting \(L\to\infty\) in the finite-projection identity therefore gives
\[
  \Var[\widehat I_N(f)]
  =\sum_{\boldsymbol\ell\in
       \Nzero^{ds}\setminus\{\boldsymbol0\}}
     \sigma_{d,\boldsymbol\ell,s}^2(f)
     \Gamma_{d,\boldsymbol\ell}(P),
\]
which is the asserted variance decomposition.
\end{proof}

For a digital \((t,m,ds)\)-net, the proof of \cite[Lemma~8]{Dick2011} remains valid when the active-coordinate set is allowed to be an arbitrary subset of \(\{1,\ldots,ds\}\). In particular, let \(M=m-t\), and let
\[
u=\supp(\boldsymbol\ell)
  :=\{j\in\{1,\ldots,ds\}:\ell_j>0\}.
\]
Then, for every
\(\boldsymbol\ell\in\Nzero^{ds}\setminus\{\boldsymbol0\}\),
\begin{equation}\label{eq:gain-piecewise}
 \Gamma_{d,\boldsymbol\ell}(P)\le
 \begin{cases}
  0,
     &\abs{\boldsymbol\ell}_1\le M,\\
  b^{\abs u-\abs{\boldsymbol\ell}_1},
     &M<\abs{\boldsymbol\ell}_1\le M+\abs u,\\
  b^{-M},
     &\abs{\boldsymbol\ell}_1>M+\abs u.
 \end{cases}
\end{equation}
Consequently,
\begin{equation}\label{eq:gain-uniform}
 \Gamma_{d,\boldsymbol\ell}(P)
 \le b^{ds}b^{-M}
       \one_{\{\abs{\boldsymbol\ell}_1>M\}}.
\end{equation}

\section{The finite-difference variation and the Sobolev norm}

\subsection{The asserted equality is false}

The statement on pp.~1386--1387 that the finite-difference variation coincides
with the displayed Sobolev norm is withdrawn.

\begin{example}[A smooth counterexample]\label{ex:variation-counterexample}
Take $b=2$, $s=d=\alpha=1$, and $f(x)=x^2/2$.  In the definition of the
one-dimensional Vitali term, use the partition
$[0,1/2)\cup[1/2,1)$.  On the first cell, with increment $z=1/2$,
\[
 \frac{f(t+1/2)-f(t)}{1/2}=t+\frac14,
\]
whose supremum is $3/4$.  On the second cell the increment $z=-1/2$ gives the
same supremum.  Hence the square of the finite-difference Vitali term is at
least
\[
 2\cdot\frac12\left(\frac34\right)^2=\frac9{16}.
\]
The square of the claimed derivative expression is
\[
 \int_0^1\abs{f'(x)}^2\dd x=\int_0^1x^2\dd x=\frac13.
\]
Thus the two quantities are not equal even for a polynomial.
\end{example}

The fine uniform partitions do converge to the derivative integral, but the
variation is a supremum over all admissible partitions, including coarse ones.
A limit over a refining sequence cannot therefore be substituted for that
supremum.

\subsection{A second obstruction in Appendix B}

There is also an independent admissibility problem in the last step of
Appendix~B.  The variation $V_\alpha$ only permits, within one output
coordinate, increments whose base-$b$ positions have the form
\[
 \alpha(L_r-1)+r,\qquad r=1,\ldots,\alpha.
\]
The increments produced by an order-$d$ Walsh block have positions
\[
 h+d(\ell_h-1),\qquad h=1,\ldots,d.
\]
After retaining $\min(\alpha,|K_i|)$ increments, these positions need not have
the residue pattern required by the definition of $V_\alpha$.  For example,
with $d=3$, $\alpha=2$, $s=1$, $K=\{1,3\}$ and
$\ell_1=\ell_3=1$, the two positions are $1$ and $3$, both odd.  An admissible
order-two variation cell requires one odd and one even position.  Thus the
last displayed inequality on p.~1397 does not follow from the definition of
$V_\alpha$.

This observation does not prove that every possible estimate in terms of
$V_\alpha$ is false, but it shows that Lemma~9 and Theorem~10 are not established
for that finite-difference variation.  The variation version is therefore
withdrawn.  The direct Sobolev estimate below is independent of both defects.

\subsection{The unanchored mixed Sobolev norm}

For $u\subseteq\{1,\ldots,s\}$ define the marginal
\[
 f_u(x_u):=\int_{[0,1]^{s-\abs u}}f(x)\dd x_{-u}.
\]
For a smooth function, define
\begin{equation}\label{eq:sobolev-norm}
 \norm{f}_{\mathcal H_{s,\alpha}}^2
 :=\abs{I(f)}^2
 +\sum_{\emptyset\ne u\subseteq\{1,\ldots,s\}}
   \ \sum_{\boldsymbol\nu_u\in\{1,\ldots,\alpha\}^{u}}
   \int_{[0,1]^{\abs u}}
   \abs{\partial^{\boldsymbol\nu_u}f_u(x_u)}^2\dd x_u.
\end{equation}
This is the norm displayed after the definition of $V_\alpha(f)$ in the
original paper.  The notation $\mathcal H_{s,\alpha}$ is introduced here to
avoid identifying it with the finite-difference variation.

The norm controls $L_2$.  Indeed, for the ANOVA decomposition
$f=\sum_ug_u$, orthogonality gives
$\norm f_2^2=\sum_u\norm{g_u}_2^2$.  For $u\ne\emptyset$, $g_u$ has zero
mean in every active coordinate, and iterated one-dimensional Poincare
inequalities give
\[
 \norm{g_u}_2\le\norm{\partial^{\boldsymbol1_u}g_u}_2
 =\norm{\partial^{\boldsymbol1_u}f_u}_2.
\]
Hence
\begin{equation}\label{eq:H-embeds-L2}
 \norm f_{L_2([0,1]^s)}\le\norm f_{\mathcal H_{s,\alpha}}.
\end{equation}
The completion of smooth functions in \eqref{eq:sobolev-norm} is therefore
canonically embedded in $L_2$.  All coefficient estimates below extend to that
completion by density.

\section{Corrected decay of the Walsh blocks}

Let $\boldsymbol\ell=(\ell_1,\ldots,\ell_{ds})\in\Nzero^{ds}$ and define
\[
 K:=\set{j:\ell_j>0},\qquad
 K_i:=K\cap\set{(i-1)d+1,\ldots,id}.
\]
Put
\[
 v:=\set{i:K_i\ne\emptyset},\qquad
 r_i:=\min(\alpha,\abs{K_i})\quad(i\in v),
\]
and
\begin{equation}\label{eq:q-def}
 q(\boldsymbol\ell):=\sum_{i\in v}(\abs{K_i}-r_i).
\end{equation}
For $j=(i-1)d+h\in K_i$, $1\le h\le d$, define
\begin{equation}\label{eq:eta-j}
 \eta'_j:=(b-1)b^{-h-(\ell_j-1)d}
          =(b-1)b^{d-h}b^{-d\ell_j}.
\end{equation}
For each $i\in v$, arrange all $\abs{K_i}$ numbers in nondecreasing order,
with repetitions retained:
\[
 \eta_{i,1}\le\cdots\le\eta_{i,\abs{K_i}},
 \qquad
 \set{\eta_{i,1},\ldots,\eta_{i,\abs{K_i}}}
 =\set{\eta'_j:j\in K_i}
\]
as multisets, and put
\begin{equation}\label{eq:eta-block}
 \eta(\boldsymbol\ell)
 :=\prod_{i\in v}\prod_{h=1}^{r_i}\eta_{i,h}.
\end{equation}
This corrects the impossible equality in the printed Lemma~9, where an
$r_i$-element set was identified with an $\abs{K_i}$-element set.

\begin{lemma}[Order statistics]\label{lem:order-statistics}
Let $0\le a_j\le c_j$ for $1\le j\le n$, and write
$a_{(1)}\le\cdots\le a_{(n)}$ and
$c_{(1)}\le\cdots\le c_{(n)}$ for the corresponding order statistics, with
repetitions.  Then $a_{(h)}\le c_{(h)}$ for every $h$, and therefore
\[
 \prod_{h=1}^r a_{(h)}\le\prod_{h=1}^r c_{(h)}
 \qquad(0\le r\le n).
\]
\end{lemma}

\begin{proof}
If $a_{(h)}>c_{(h)}$, at least $h$ indices satisfy
$c_j\le c_{(h)}<a_{(h)}$.  For those indices,
$a_j\le c_j<a_{(h)}$, contradicting the definition of the $h$th order
statistic.  Multiplication gives the second assertion.
\end{proof}

\begin{lemma}[Direct Sobolev replacement for Lemma 9]\label{lem:block-decay}
If $f\in\mathcal H_{s,\alpha}$, then
\begin{equation}\label{eq:sobolev-block-local}
 \sigma_{d,\boldsymbol\ell,s}(f)
 \le 2^{q(\boldsymbol\ell)}b^{\abs K}
      \eta(\boldsymbol\ell)
      \norm{\partial^{\boldsymbol r}f_v}_{L_2([0,1]^{\abs v})},
\end{equation}
where $\boldsymbol r=(r_i)_{i\in v}$.  Consequently,
\begin{equation}\label{eq:sobolev-block}
 \sigma_{d,\boldsymbol\ell,s}(f)
 \le 2^{s\max(d-\alpha,0)}b^{ds}
      \eta(\boldsymbol\ell)\norm f_{\mathcal H_{s,\alpha}}.
\end{equation}
\end{lemma}

\begin{proof}
Write
\[
 L:=|\boldsymbol\ell|_1
   =\sum_{j=1}^{ds}\ell_j,
 \qquad
 K:=\{j\in\{1,\ldots,ds\}:\ell_j>0\}.
\]
For each output coordinate $i\in\{1,\ldots,s\}$, define
\[
 K_i
 :=
 K\cap\{(i-1)d+1,\ldots,id\},
 \qquad
 v:=\{i:K_i\neq\emptyset\},
\]
and put
\[
 r_i:=\min(\alpha,|K_i|),
 \qquad
 q(\boldsymbol\ell)
 :=
 \sum_{i\in v}(|K_i|-r_i).
\]
Thus $r_i$ is the number of increments which will be retained in
output coordinate $i$, while $q(\boldsymbol\ell)$ is the total
number of discarded increments.

We first derive the finite Walsh-block identity which is the starting
point of the argument.

For
$\boldsymbol m=(m_1,\ldots,m_{ds})\in\Nzero^{ds}$, define
\[
 \mathcal A_{\boldsymbol m}
 :=
 \prod_{j=1}^{ds}\{0,\ldots,b^{m_j}-1\}.
\]
In particular,
\[
 \mathcal A_{\boldsymbol\ell}
 =
 \prod_{j=1}^{ds}\{0,\ldots,b^{\ell_j}-1\},
 \qquad
 \mathcal A_1:=\{0,\ldots,b-1\}^{ds}.
\]
For $u\subseteq K$, let $\boldsymbol1_u$ denote its indicator
vector.  Since, for $j\in K$,
\[
 \mathbf 1_{\{b^{\ell_j-1}\le k_j<b^{\ell_j}\}}
 =
 \mathbf 1_{\{0\le k_j<b^{\ell_j}\}}
 -
 \mathbf 1_{\{0\le k_j<b^{\ell_j-1}\}},
\]
multiplication over the active coordinates gives
\begin{equation}\label{eq:block-indicator-detailed}
 \mathbf1_{B_{\boldsymbol\ell}}(\boldsymbol k)
 =
 \sum_{u\subseteq K}
 (-1)^{|u|}
 \mathbf1_{\mathcal A_{\boldsymbol\ell-\boldsymbol1_u}}
 (\boldsymbol k).
\end{equation}
Consequently,
\begin{equation}\label{eq:block-inclusion-exclusion-detailed}
 \sigma_{d,\boldsymbol\ell,s}^2(f)
 =
 \sum_{u\subseteq K}(-1)^{|u|}
 \sum_{\boldsymbol k\in
       \mathcal A_{\boldsymbol\ell-\boldsymbol1_u}}
 \left|
   \widehat f(\E_d(\boldsymbol k))
 \right|^2.
\end{equation}

Equivalently, define the finite Walsh polynomial
\[
 G_{\boldsymbol\ell}(x)
 :=
 \sum_{\boldsymbol k\in B_{\boldsymbol\ell}}
 \widehat f(\E_d(\boldsymbol k))
 \wal_{\E_d(\boldsymbol k)}(x).
\]
By orthogonality of the Walsh system,
\begin{equation}\label{eq:block-parseval-detailed}
 \sigma_{d,\boldsymbol\ell,s}^2(f)
 =
 \int_{[0,1]^s}|G_{\boldsymbol\ell}(x)|^2\,\dd x.
\end{equation}
Using \eqref{eq:block-indicator-detailed}, we may also write
\begin{equation}\label{eq:G-inclusion-exclusion-detailed}
 G_{\boldsymbol\ell}(x)
 =
 \sum_{u\subseteq K}(-1)^{|u|}
 \sum_{\boldsymbol k\in
       \mathcal A_{\boldsymbol\ell-\boldsymbol1_u}}
 \widehat f(\E_d(\boldsymbol k))
 \wal_{\E_d(\boldsymbol k)}(x).
\end{equation}

All frequencies occurring in this expression vanish in the output
coordinates outside $v$.  After integrating out those inactive
coordinates, the corresponding Walsh coefficients may therefore be
computed from $f_v$.  We work from now on on
$[0,1]^{|v|}$.

For $a\in\mathcal A_{\boldsymbol\ell}$, define the fine interlaced
cell
\[
 E_a
 :=
 \operatorname{pr}_v
 \D_d\bigl(
 [ab^{-\boldsymbol\ell},
  (a+1)b^{-\boldsymbol\ell})
 \bigr).
\]
The cells $E_a$, $a\in\mathcal A_{\boldsymbol\ell}$, are pairwise
disjoint up to their boundaries and partition $[0,1]^{|v|}$.  They
fix exactly
\[
 L=\sum_{j\in K}\ell_j
\]
base-$b$ digits, and hence
\begin{equation}\label{eq:Ea-volume-detailed}
 \Vol(E_a)=b^{-L}.
\end{equation}

For $u\subseteq K$, put
\[
 \boldsymbol m_u:=\boldsymbol\ell-\boldsymbol1_u
\]
and define $\pi_u(a)\in\mathcal A_{\boldsymbol m_u}$ by
\[
 [\pi_u(a)]_j
 :=
 \begin{cases}
  \lfloor a_j/b\rfloor, & j\in u,\\
  a_j, & j\notin u.
 \end{cases}
\]
Let
\[
 E_{u,a}
 :=
 \operatorname{pr}_v
 \D_d\bigl(
 [\pi_u(a)b^{-\boldsymbol m_u},
  (\pi_u(a)+1)b^{-\boldsymbol m_u})
 \bigr).
\]
Thus $E_{u,a}$ is obtained from $E_a$ by forgetting the last
prescribed digit in each input coordinate $j\in u$.

We shall use the finite Walsh-kernel identity
\begin{equation}\label{eq:finite-walsh-kernel-detailed}
 \sum_{\boldsymbol k\in\mathcal A_{\boldsymbol m}}
 \wal_{\E_d(\boldsymbol k)}(x)
 \overline{\wal_{\E_d(\boldsymbol k)}(y)}
 =
 b^{|\boldsymbol m|_1}
 \mathbf1_{\{x,y\text{ agree in all digits prescribed by }
                  \boldsymbol m\}}.
\end{equation}
This follows coordinatewise from the one-dimensional identity
\[
 \sum_{k=0}^{b^m-1}
 \wal_k(x)\overline{\wal_k(y)}
 =
 b^m
 \mathbf1_{\{\lfloor b^m x\rfloor=\lfloor b^m y\rfloor\}},
\]
together with the fact that $\E_d$ interlaces the corresponding
frequency digits in the same way that $\D_d$ interlaces the point
digits.

For $x\in E_a$, applying
\eqref{eq:finite-walsh-kernel-detailed} gives
\begin{align}
 &\sum_{\boldsymbol k\in
       \mathcal A_{\boldsymbol\ell-\boldsymbol1_u}}
 \widehat f(\E_d(\boldsymbol k))
 \wal_{\E_d(\boldsymbol k)}(x)
 \notag\\
 &\qquad
 =
 b^{L-|u|}
 \int_{E_{u,a}} f_v(t)\,\dd t.
 \label{eq:lower-projection-cell-detailed}
\end{align}
Hence $G_{\boldsymbol\ell}$ is constant on every $E_a$, and
\begin{equation}\label{eq:G-on-cell-detailed}
 G_{\boldsymbol\ell}(x)
 =
 b^L
 \sum_{u\subseteq K}
 (-1)^{|u|}b^{-|u|}
 \int_{E_{u,a}}f_v(t)\,\dd t,
 \qquad x\in E_a.
\end{equation}
Using \eqref{eq:Ea-volume-detailed} in
\eqref{eq:block-parseval-detailed}, we obtain the exact identity
\begin{equation}\label{eq:block-cell-exact-detailed}
 \sigma_{d,\boldsymbol\ell,s}^2(f)
 =
 b^L
 \sum_{a\in\mathcal A_{\boldsymbol\ell}}
 \left|
   \sum_{u\subseteq K}
   (-1)^{|u|}b^{-|u|}
   \int_{E_{u,a}}f_v(t)\,\dd t
 \right|^2.
\end{equation}

We next resolve each parent cell $E_{u,a}$ into its fine children.
For $a\in\mathcal A_{\boldsymbol\ell}$, write
\[
 e_j(a):=b\lfloor a_j/b\rfloor,
 \qquad
 d_j(a):=a_j-e_j(a)
        =a_j-b\lfloor a_j/b\rfloor.
\]
Thus $d_j(a)\in\{0,\ldots,b-1\}$ is the last base-$b$ digit
of $a_j$.

For $u\subseteq K$ and
$\kappa_u\in\{0,\ldots,b-1\}^{u}$, let
$a^{(u,\kappa)}$ be the fine-child index defined by
\[
 a^{(u,\kappa)}_j
 :=
 \begin{cases}
  e_j(a)+\kappa_j, & j\in u,\\
  a_j, & j\notin u.
 \end{cases}
\]
Then, up to boundaries,
\[
 E_{u,a}
 =
 \bigcup_{\kappa_u\in\{0,\ldots,b-1\}^{u}}
 E_{a^{(u,\kappa)}},
\]
and therefore
\begin{equation}\label{eq:parent-children-detailed}
 \int_{E_{u,a}}f_v
 =
 \sum_{\kappa_u\in\{0,\ldots,b-1\}^{u}}
 \int_{E_{a^{(u,\kappa)}}}f_v.
\end{equation}
Since the summand depends only on the components of the replacement
digit vector in $u$, we may equivalently average over the full set
$\mathcal A_1=\{0,\ldots,b-1\}^{ds}$:
\begin{equation}\label{eq:full-digit-average-detailed}
 b^{-|u|}
 \int_{E_{u,a}}f_v
 =
 b^{-ds}
 \sum_{k\in\mathcal A_1}
 \int_{E_{a^{(u,k)}}}f_v.
\end{equation}
Indeed, each choice of the digits $k_j$, $j\in u$, occurs exactly
$b^{ds-|u|}$ times in the full sum over $k$.

We now express the child cells as translates of $E_a$.
If
\[
 j=(i-1)d+h\in K,
 \qquad 1\le h\le d,
\]
define
\begin{equation}\label{eq:z-aj-detailed}
 z_j(a,k)
 :=
 (k_j-d_j(a))b^{-h-(\ell_j-1)d}.
\end{equation}
Changing the final prescribed digit in input coordinate $j$ from
$d_j(a)$ to $k_j$ changes the corresponding digit of output
coordinate $i$, namely the digit in position
\[
 h+d(\ell_j-1),
\]
and therefore translates that output coordinate by exactly
$z_j(a,k)$.

For $u\subseteq K$, let $Z_u(a,k)\in\mathbb R^{|v|}$ have
$i$th active component
\[
 [Z_u(a,k)]_i
 :=
 \sum_{j\in u\cap K_i}z_j(a,k).
\]
Then, up to boundaries,
\[
 E_{a^{(u,k)}}=E_a+Z_u(a,k),
\]
and translation invariance of Lebesgue measure gives
\begin{equation}\label{eq:child-translate-detailed}
 \int_{E_{a^{(u,k)}}}f_v(x)\,\dd x
 =
 \int_{E_a}f_v(t+Z_u(a,k))\,\dd t.
\end{equation}

Define
\begin{equation}\label{eq:delta-explicit-detailed}
 \delta_{a,k}(t)
 :=
 \sum_{u\subseteq K}(-1)^{|u|}
 f_v\bigl(t+Z_u(a,k)\bigr).
\end{equation}
Combining
\eqref{eq:full-digit-average-detailed} and
\eqref{eq:child-translate-detailed}, we obtain
\begin{align}
 &\sum_{u\subseteq K}
 (-1)^{|u|}b^{-|u|}
 \int_{E_{u,a}}f_v
 \notag\\
 &\qquad
 =
 b^{-ds}
 \sum_{k\in\mathcal A_1}
 \int_{E_a}\delta_{a,k}(t)\,\dd t.
 \label{eq:block-difference-rearrangement-detailed}
\end{align}
Substituting this identity into
\eqref{eq:block-cell-exact-detailed} yields the exact finite
Walsh-block identity
\begin{equation}\label{eq:block-exact-detailed}
 \sigma_{d,\boldsymbol\ell,s}^2(f)
 =
 b^{L-2ds}
 \sum_{a\in\mathcal A_{\boldsymbol\ell}}
 \left|
   \sum_{k\in\mathcal A_1}
   \int_{E_a}\delta_{a,k}(t)\,\dd t
 \right|^2.
\end{equation}
Up to this point all steps are equalities; no Sobolev estimate has
yet been used.

By the triangle inequality,
\[
 \left|
   \sum_{k\in\mathcal A_1}
   \int_{E_a}\delta_{a,k}(t)\,\dd t
 \right|
 \le
 \sum_{k\in\mathcal A_1}
 \int_{E_a}|\delta_{a,k}(t)|\,\dd t.
\]
Set
\[
 A_{a,k}
 :=
 \int_{E_a}|\delta_{a,k}(t)|\,\dd t.
\]
Then
\begin{equation}\label{eq:before-max-detailed}
 \sigma_{d,\boldsymbol\ell,s}^2(f)
 \le
 b^{L-2ds}
 \sum_{a\in\mathcal A_{\boldsymbol\ell}}
 \left(
   \sum_{k\in\mathcal A_1}A_{a,k}
 \right)^2.
\end{equation}

By Cauchy--Schwarz in the finite sum over $k$, using
$|\mathcal A_1|=b^{ds}$,
\[
 \left(
   \sum_{k\in\mathcal A_1}A_{a,k}
 \right)^2
 \le
 b^{ds}\sum_{k\in\mathcal A_1}A_{a,k}^2.
\]
A second application of Cauchy--Schwarz, now on $E_a$, gives
\[
 A_{a,k}^2
 \le
 \Vol(E_a)
 \int_{E_a}|\delta_{a,k}(t)|^2\,\dd t
 =
 b^{-L}
 \int_{E_a}|\delta_{a,k}(t)|^2\,\dd t.
\]
Substituting these estimates into
\eqref{eq:before-max-detailed} gives
\begin{equation}\label{eq:block-L2-reduction-detailed}
 \sigma_{d,\boldsymbol\ell,s}^2(f)
 \le
 b^{-ds}
 \sum_{k\in\mathcal A_1}
 \sum_{a\in\mathcal A_{\boldsymbol\ell}}
 \int_{E_a}|\delta_{a,k}(t)|^2\,\dd t.
\end{equation}

We now estimate the finite differences.  Fix $a$ and $k$.
From \eqref{eq:z-aj-detailed},
\[
 |z_j(a,k)|
 \le
 (b-1)b^{-h-(\ell_j-1)d}
 =:\eta'_j,
 \qquad
 j=(i-1)d+h\in K.
\]
For a fixed output coordinate $i\in v$, arrange the numbers
$|z_j(a,k)|$, $j\in K_i$, in nondecreasing order and retain the
$r_i$ smallest.  Finite-difference operators in the same coordinate
commute, so the order of the increments is irrelevant.

Likewise arrange the numbers $\eta'_j$, $j\in K_i$, in
nondecreasing order.  Since
\[
 |z_j(a,k)|\le\eta'_j,
\]
Lemma~\ref{lem:order-statistics} shows that the $h$th smallest actual
increment is no larger than the $h$th smallest upper bound.  Hence
the product of the $r_i$ retained absolute increments is at most
the product of the $r_i$ smallest $\eta'_j$.  Multiplication over
$i\in v$, together with the definition \eqref{eq:eta-block}, gives
\begin{equation}\label{eq:retained-product-detailed}
 \prod_{i\in v}
 \prod_{\substack{j\in K_i\\ j\ {\rm retained}}}
 |z_j(a,k)|
 \le
 \eta(\boldsymbol\ell).
\end{equation}

We next reduce the order of the finite difference.  Write
\[
 T_yg(t):=g(t+y).
\]
Up to the harmless global sign determined by the convention in
\eqref{eq:delta-explicit-detailed}, the contribution from output
coordinate $i$ is
\[
 \prod_{j\in K_i}(I-T_{z_j e_i}).
\]
Let $R_i\subseteq K_i$ be the indices of the $r_i$ retained
increments and let
\[
 Q_i:=K_i\setminus R_i.
\]
Then
\begin{align}
 \prod_{j\in K_i}(I-T_{z_j e_i})
 &=
 \left(
   \prod_{j\in Q_i}(I-T_{z_j e_i})
 \right)
 \left(
   \prod_{j\in R_i}(I-T_{z_j e_i})
 \right)
 \notag\\
 &=
 \sum_{U_i\subseteq Q_i}
 (-1)^{|U_i|}
 T_{\sum_{j\in U_i}z_j e_i}
 \prod_{j\in R_i}(I-T_{z_j e_i}).
 \label{eq:difference-reduction-detailed}
\end{align}
Thus a difference of order $|K_i|$ is an alternating sum of
$2^{|K_i|-r_i}$ translates of differences of order $r_i$.
Multiplication over all active output coordinates shows that
$\delta_{a,k}$ is an alternating sum of exactly
\[
 \prod_{i\in v}2^{|K_i|-r_i}
 =
 2^{q(\boldsymbol\ell)}
\]
translates of mixed finite differences of order
$\boldsymbol r=(r_i)_{i\in v}$.

For a signed increment $z$, write
\[
 I(z):=[\min(0,z),\max(0,z)].
\]
Repeated use of the fundamental theorem of calculus gives, up to an
irrelevant sign,
\begin{align}
 \Delta_{\boldsymbol r}g(t;\boldsymbol z)
 &=
 \int_{I(z_{1,1})}\!\!\cdots
 \int_{I(z_{i,h})}\!\!\cdots
 \partial^{\boldsymbol r}g
 \left(
   t+\sum_{i,h}\theta_{i,h}e_i
 \right)
 \prod_{i,h}\dd\theta_{i,h}.
 \label{eq:FTC-mixed-detailed}
\end{align}
Consequently, if $S\subseteq[0,1]^{|v|}$ and the translated points
which occur in the integral remain in $[0,1]^{|v|}$, Minkowski's
integral inequality gives
\begin{align}
 \left\|
 \Delta_{\boldsymbol r}g(\,\cdot\,;\boldsymbol z)
 \right\|_{L_2(S)}
 &\le
 \int_{I(z_{1,1})}\!\!\cdots
 \int_{I(z_{i,h})}\!\!\cdots
 \left\|
 \partial^{\boldsymbol r}g
 \left(
   \,\cdot\,
   +\sum_{i,h}\theta_{i,h}e_i
 \right)
 \right\|_{L_2(S)}
 \prod_{i,h}\dd\theta_{i,h}
 \notag\\
 &\le
 \left(\prod_{i,h}|z_{i,h}|\right)
 \left\|
   \partial^{\boldsymbol r}g
 \right\|_{L_2([0,1]^{|v|})}.
 \label{eq:finite-difference-L2-detailed}
\end{align}

It remains to organize the sum over the fine cells $E_a$.  The
increments in \eqref{eq:z-aj-detailed} depend on $a$ only through
the last active digits
\[
 \rho(a):=(d_j(a))_{j\in K}
 \in\{0,\ldots,b-1\}^{K}.
\]
For
\[
 \rho\in\{0,\ldots,b-1\}^{K},
\]
define
\[
 S_\rho
 :=
 \bigcup_{\substack{
 a\in\mathcal A_{\boldsymbol\ell}\\
 d_j(a)=\rho_j\ {\rm for\ all}\ j\in K}}
 E_a.
\]
The union is disjoint up to boundaries.  If $\rho(a)=\rho$, then,
for fixed $k$, all increments $z_j(a,k)$ are the same; denote them
by $z_j(\rho,k)$.  Hence the finite-difference expression
$\delta_{a,k}$ is the same function of its argument for every
$a$ in this group.  Denote it by $\delta_{\rho,k}$.  Then
\begin{equation}\label{eq:group-last-digits-detailed}
 \sum_{\substack{
 a\in\mathcal A_{\boldsymbol\ell}\\
 \rho(a)=\rho}}
 \int_{E_a}|\delta_{a,k}(t)|^2\,\dd t
 =
 \int_{S_\rho}|\delta_{\rho,k}(t)|^2\,\dd t.
\end{equation}

For $t\in S_\rho$, every vertex which occurs in
$\delta_{\rho,k}(t)$ belongs to $[0,1]^{|v|}$.  Indeed, adding
$z_j(\rho,k)$ changes the relevant prescribed interlaced digit
from $\rho_j$ to $k_j$ and leaves all other prescribed digits
unchanged.  The boxes parametrized by the variables
$\theta_{i,h}$ in \eqref{eq:FTC-mixed-detailed} are contained in
the convex hull of the corresponding finite-difference vertices.
Since the cube $[0,1]^{|v|}$ is convex, these boxes remain in the
cube.

We may therefore apply
\eqref{eq:finite-difference-L2-detailed} to each of the
$2^{q(\boldsymbol\ell)}$ reduced differences.  By the triangle
inequality in $L_2$, followed by
\eqref{eq:retained-product-detailed},
\begin{equation}\label{eq:delta-L2-detailed}
 \|\delta_{\rho,k}\|_{L_2(S_\rho)}
 \le
 2^{q(\boldsymbol\ell)}
 \eta(\boldsymbol\ell)
 \left\|
   \partial^{\boldsymbol r}f_v
 \right\|_{L_2([0,1]^{|v|})}.
\end{equation}
Squaring and using
\eqref{eq:group-last-digits-detailed}, we obtain
\[
 \sum_{\substack{
 a\in\mathcal A_{\boldsymbol\ell}\\
 \rho(a)=\rho}}
 \int_{E_a}|\delta_{a,k}(t)|^2\,\dd t
 \le
 4^{q(\boldsymbol\ell)}
 \eta(\boldsymbol\ell)^2
 \left\|
   \partial^{\boldsymbol r}f_v
 \right\|_2^2.
\]

There are exactly $b^{|K|}$ possible vectors
$\rho\in\{0,\ldots,b-1\}^{K}$.  Therefore, for every fixed
$k\in\mathcal A_1$,
\begin{equation}\label{eq:sum-a-bound-detailed}
 \sum_{a\in\mathcal A_{\boldsymbol\ell}}
 \int_{E_a}|\delta_{a,k}(t)|^2\,\dd t
 \le
 b^{|K|}
 4^{q(\boldsymbol\ell)}
 \eta(\boldsymbol\ell)^2
 \left\|
   \partial^{\boldsymbol r}f_v
 \right\|_2^2.
\end{equation}
Substituting \eqref{eq:sum-a-bound-detailed} into
\eqref{eq:block-L2-reduction-detailed} and using
$|\mathcal A_1|=b^{ds}$, we obtain
\begin{align*}
 \sigma_{d,\boldsymbol\ell,s}^2(f)
 &\le
 b^{-ds}b^{ds}b^{|K|}
 4^{q(\boldsymbol\ell)}
 \eta(\boldsymbol\ell)^2
 \left\|
   \partial^{\boldsymbol r}f_v
 \right\|_2^2\\
 &=
 4^{q(\boldsymbol\ell)}b^{|K|}
 \eta(\boldsymbol\ell)^2
 \left\|
   \partial^{\boldsymbol r}f_v
 \right\|_2^2.
\end{align*}
Taking square roots yields
\begin{equation}\label{eq:sobolev-block-local-strong-detailed}
 \sigma_{d,\boldsymbol\ell,s}(f)
 \le
 2^{q(\boldsymbol\ell)}b^{|K|/2}
 \eta(\boldsymbol\ell)
 \left\|
   \partial^{\boldsymbol r}f_v
 \right\|_2.
\end{equation}
Since $b^{|K|/2}\le b^{|K|}$, this implies
\[
 \sigma_{d,\boldsymbol\ell,s}(f)
 \le
 2^{q(\boldsymbol\ell)}b^{|K|}
 \eta(\boldsymbol\ell)
 \left\|
   \partial^{\boldsymbol r}f_v
 \right\|_{L_2([0,1]^{|v|})},
\]
which proves \eqref{eq:sobolev-block-local}.

Finally,
\[
 |K_i|-r_i
 =
 |K_i|-\min(\alpha,|K_i|)
 =
 \max(|K_i|-\alpha,0)
 \le
 \max(d-\alpha,0),
\]
because $|K_i|\le d$.  Summing over at most $s$ output coordinates
gives
\[
 q(\boldsymbol\ell)
 \le
 s\max(d-\alpha,0).
\]
Moreover,
\[
 |K|\le ds.
\]
By the definition of the Sobolev norm in
\eqref{eq:sobolev-norm},
\[
 \left\|
   \partial^{\boldsymbol r}f_v
 \right\|_{L_2([0,1]^{|v|})}
 \le
 \|f\|_{\mathcal H_{s,\alpha}}.
\]
Consequently,
\[
 \sigma_{d,\boldsymbol\ell,s}(f)
 \le
 2^{s\max(d-\alpha,0)}
 b^{ds}
 \eta(\boldsymbol\ell)
 \|f\|_{\mathcal H_{s,\alpha}},
\]
which is \eqref{eq:sobolev-block}.
\end{proof}

\begin{remark}
The $b^{\abs K}$ factor in \eqref{eq:sobolev-block-local} is deliberately kept
as a simple, safe constant.  It has no effect on the rate in $N$.  It also makes the Sobolev repair independent of the two defects in the
finite-difference variation argument.
Sharper constants may be obtained by retaining the averaging over the inner
digit indices, but are not needed for the result of the paper.
\end{remark}

\section{Corrected main theorem}

Throughout this section put
\begin{equation}\label{eq:kappa}
 \kappa(d,s):=ds-1.
\end{equation}

\begin{lemma}[Exponential polynomial tail]\label{lem:tail}
Let $0<q<1$ and $r\in\Nzero$.  There is a finite constant $C_{q,r}$ such that,
for every $L\in\Nzero$,
\[
 \sum_{p=L+1}^{\infty}q^p(p+2)^r
 \le C_{q,r}\,q^L(L+2)^r.
\]
\end{lemma}

\begin{proof}
Write $p=L+n$.  Since $L+n+2\le(L+2)(n+1)$,
\[
 \sum_{n=1}^{\infty}q^{L+n}(L+n+2)^r
 \le q^L(L+2)^r\sum_{n=1}^{\infty}q^n(n+1)^r.
\]
The last series is finite.
\end{proof}

\begin{theorem}[Corrected Theorem 10]\label{thm:main}
Let $P$ be a digital $(t,m,ds)$-net in base $b$, and let
$\widehat I_N$ be its order-$d$ nested-uniformly scrambled and interlaced
estimator.  There is a constant $C_{b,s,\alpha,d}>0$, independent of $m$ and
$f$, such that, for every $f\in\mathcal H_{s,\alpha}$,
\begin{equation}\label{eq:main-sobolev}
 \Var[\widehat I_N(f)]
 \le C_{b,s,\alpha,d}\,\norm f_{\mathcal H_{s,\alpha}}^2
 b^{-(2\rho+1)M}(M+2)^{ds-1}.
\end{equation}
\end{theorem}

\begin{proof}
Use Lemma~\ref{lem:variance-decomposition}, the gain bound
\eqref{eq:gain-uniform}, and Lemma~\ref{lem:block-decay}.  All
factors depending only on $b,s,\alpha,d$ are absorbed into
$C_{b,s,\alpha,d}$.

Suppose first that $d\le\alpha$.  Then $r_i=\abs{K_i}$ and
$q(\boldsymbol\ell)=0$.  From \eqref{eq:eta-j},
\[
 \eta(\boldsymbol\ell)
 \le A_{b,d,s}\,b^{-d\abs{\boldsymbol\ell}_1}.
\]
The number of vectors in $\Nzero^{ds}$ with coordinate sum $p$ is
$\binom{p+ds-1}{ds-1}$.  Hence Lemma~\ref{lem:tail} gives
\[
 \sum_{\abs{\boldsymbol\ell}_1>M}
 \eta(\boldsymbol\ell)^2
 \le C_{b,d,s}\,b^{-2dM}(M+2)^{ds-1}.
\]
Multiplication by the factor $b^{-M}$ from \eqref{eq:gain-uniform} proves the
claim in this case.

Now suppose that $d>\alpha$.  Within each group of $d$ coordinates, arrange the
levels in nonincreasing order,
$\ell_{i,1}\ge\cdots\ge\ell_{i,d}$.  There are at most $(d!)^s$ reorderings and
\[
 \eta(\boldsymbol\ell)
 \le A_{b,d,\alpha,s}
 b^{-d\sum_{i=1}^s\sum_{j=1}^{\alpha}\ell_{i,j}}.
\]
If $\sum_{i,j}\ell_{i,j}>M$, then, since the mean of the largest $\alpha$
entries is at least the mean of all $d$ entries,
\[
 \sum_{i=1}^s\sum_{j=1}^{\alpha}\ell_{i,j}>\frac{\alpha M}{d}.
\]
Put $k_i=\sum_{j=1}^{\alpha}\ell_{i,j}$ and $p_i=dk_i$.  For fixed $k_i$, the
number of choices of the first $\alpha$ entries is at most
$\binom{k_i+\alpha-1}{\alpha-1}$, while each of the remaining $d-\alpha$
entries is at most $\floor{k_i/\alpha}$.  After relaxing the condition that $p_i$ is divisible by $d$, the required
sum is bounded by
\[
 C\sum_{\substack{p_1,\ldots,p_s\ge0\\p_1+\cdots+p_s>\alpha M}}
 b^{-2(p_1+\cdots+p_s)}
 \prod_{i=1}^s(p_i+2)^{d-1}.
\]
Grouping by $p=p_1+\cdots+p_s$ and using
$\binom{p+s-1}{s-1}\le(p+2)^{s-1}$ gives
\[
 C\sum_{p=\alpha M+1}^{\infty}
 b^{-2p}(p+2)^{sd-1}.
\]
Lemma~\ref{lem:tail} bounds this by
\[
 Cb^{-2\alpha M}(M+2)^{ds-1}.
\]
Again multiply by $b^{-M}$ from \eqref{eq:gain-uniform}.  This proves \eqref{eq:main-sobolev}.
\end{proof}

\begin{corollary}[Root mean square rate]\label{cor:rmse}
Under the assumptions of Theorem~\ref{thm:main},
\[
 \left(\mathbb E\abs{\widehat I_N(f)-I(f)}^2\right)^{1/2}
 \le C_{b,s,\alpha,d}^{1/2}\norm f_{\mathcal H_{s,\alpha}}
 b^{-(\rho+1/2)(m-t)}(m-t+2)^{(ds-1)/2}.
\]
For a family with fixed $b,s,d,\alpha,t$,
\[
 \operatorname{RMSE}
 =O\!\left(N^{-\rho-1/2}(\log N)^{(ds-1)/2}\right)
 =O\!\left(N^{-\rho-1/2+\eps}\right)
\]
for every $\eps>0$.  In particular, if $d\ge\alpha$, the algebraic rate is
$N^{-\alpha-1/2+\eps}$.
\end{corollary}

\section{Line-by-line corrections}\label{sec:line-corrections}

The following corrections are in addition to the replacement statements above.
Page and line descriptions refer to the journal version.

\begin{enumerate}[leftmargin=2.5em]
\item Page 1380, Proposition~4(1): ``$q$-adic rational'' should read
      ``$b$-adic rational.''
\item Page 1380: the domain of $\E_d$ is $\Nzero^d$, and in $s$ dimensions it is
      $\Nzero^{ds}$.
\item Pages 1381--1382: replace the use of $\D_d^{-1}$ by
      \eqref{eq:correct-order-d-scramble}; replace Lemma~6 by
      Proposition~\ref{prop:covariance}.
\item Page 1382, equation~(3.2): omit the zero frequency, as in
      \eqref{eq:parseval-var}.
\item Page 1387, Lemma~9: sort all $\abs{K_i}$ values and take the first
      $r_i=\min(\alpha,\abs{K_i})$, as in \eqref{eq:eta-block}.
\item Page 1388, Theorem~10: the right-hand side must be homogeneous of
      degree two in $f$.  The theorem for the original finite-difference
      variation is withdrawn and replaced by Theorem~\ref{thm:main}.  Write
      $C_{b,s,\alpha,d}$; the use of $M+2$ also covers $m=t$.
\item Page 1390, Section~4: replace the displayed logarithmic exponent by $(ds-1)/2$ as in Corollary~\ref{cor:rmse}.  The statement
      $O(N^{-\min(\alpha,d)-1/2+\eps})$ is unchanged.
\item Page 1391, proof of Lemma~12: in the case $s=1$, the vector is
      $\boldsymbol\nu=(\nu_1,\ldots,\nu_d)$, the maxima run over
      $1\le i\le d$, and $\abs{\boldsymbol\nu}_1=\sum_{i=1}^d\nu_i$.
\item Page 1392, Appendix~B: $q$ has $ds$ components, not $\alpha s$
      components; all products defining the input cells run to $ds$.
\item Page 1393: $\sigma_{d,\boldsymbol\ell,s,r}$ is
      $\sigma_{d,\boldsymbol\ell,s}$, and the Plancherel integral is over
      $[0,1]^s$.
\item Page 1394: in the definition of $\delta_k$, the zero vector is indexed by
      $\{1,\ldots,ds\}\setminus u$, not by $\{1,\ldots,\alpha s\}\setminus u$.
\item Page 1395, Lemma~13: first form the full list of $\abs{K_i}$ increments;
      then retain $r_i$ of them.  
\item Page 1396: The increments $z_{v}$ are ordered by nondecreasing
      absolute value, not by signed value.
\item Page 1396: for $n=\abs K>\alpha$ the correct identity is
\[
 \Delta_n(t;z_1,\ldots,z_n)
 =\sum_{u\subseteq\{\alpha+1,\ldots,n\}}
   (-1)^{\abs u}
   \Delta_\alpha\!\left(t+\sum_{j\in u}z_j;z_1,\ldots,z_\alpha\right),
\]
      and therefore the factor is $2^{n-\alpha}$, not $2^{\alpha-n}$.
\item Page 1397: after squaring Lemma~13, the factor is
      $4^{q(\boldsymbol\ell)}$ (bounded by $4^{s(d-\alpha)}$ when
      $d\ge\alpha$), not $2^{s(d-\alpha)}$.  The denominator contains the
      product of the retained $r_i$ increments only.  The subscript
      $\sigma_{\alpha,\boldsymbol\ell,s}$ is
      $\sigma_{d,\boldsymbol\ell,s}$.
\item Page 1397, final paragraph: when $d<\alpha$, apply the same proof directly
      with $r_i=\abs{K_i}$, or equivalently use the smoothness parameter $d$ in
      the preceding argument.  The phrase ``with $d=\alpha$'' is incorrect.
\end{enumerate}

\end{document}